\documentclass[oneside,10pt,a4paper]{book}
\usepackage[cp1251]{inputenc}
\usepackage[russian]{babel}
\usepackage{cite}
\usepackage{amscd}
\usepackage[unicode]{hyperref}
\usepackage{amssymb,amsmath,amsfonts,graphicx,makeidx,verbatim}
\makeindex
\usepackage{srcltx,fixltx2e}
\usepackage{multicol}
\begin{document}

УДК 512.554

\begin{center} {\bf Тернарные дифференцирования йордановых cупералгебр}
\end{center}

\begin{center} {\bf А.И. Шестаков}
\end{center}

{\bf Аннотация:} \emph{В работе описываются тернарные и обобщенные
дифференцирования конечномерных полупростых йордановых супералгебр
над алгебраически замкнутым полем характеристики ноль, а также простых йордановых
супералгебр с полупростой четной частью над алгебраически замкнутым полем
произвольной харакеристики не равной 2. }

{\bf Ключевые слова:} \emph{супералгебра, йорданова алгебра,
обобщенное дифференцирование, тернарное дифференцирование. }
\vspace{6pt}

\footnotetext{Работа поддержана грантом РФФИ 11-01-00938-а}

Понятие тернарного дифференцирования восходит к так называемоему
\emph{принципу тройственности} Джекобсона \cite{Schafer66},
который утверждает, что каждый элемент ортогональной алгебры Ли
кососимметрических относительно композиционной формы операторов
алгебры Кэли-Диксона является главной компонентой тернарного
дифференцирования. Непосредственно понятие тернарного
дифференцирования было введено К. Хименес-Хестал и Х.
Перес-Искьердо в \cite{JG3}. Ими описаны тернарные
дифференцирования обобщенных алгебр Кэли-Диксона над полем
характеристики не равной 2,3. В работах \cite{JG8,P-Izq9}
тернарные дифференцирования использовались для изучения некоторых
неассоцитивных унитальных алгебр.

Тернарное дифференцирование естественным образом обобщает обычные
дифференцирования, а также более широкий класс
$\delta$-дифференцирований, введенных В.Т.Филипповым в
\cite{Fill98}. В дальнейшем понятие $\delta$-дифференцирования
было развито и исследовано для многих классов алгебр и супералгебр
в работах \cite{Fill99}-\cite{kay_de5}.

Тернарные дифференцирования, в свою очередь, тесно связаны с
обобщенными дифференцированиями из \cite{LL}. Для ассоциативных
алгебр обобщенные дифференцирования были введены \mbox{М.
Брешаром} \cite{Bresar91} и \mbox{Х. Коматсу}, \mbox{А.
Накаджимой} \cite{komatsu_nakajima03}. Также в
\cite{komatsu_nakajima04} изучались ассоциативные алгебры, у
которых множество значений обобщенного дифференцирования состоит
из обратимых элементов и нуля. Для алгебр Ли обобщенные
дифференцирования изучались Дж. Леджером, Е. Лаксом в \cite{{LL}},
и для супералгебр Ли --- \mbox{Р. Жангом}, \mbox{Ю. Жангом} в
\cite{ZhangRY}. Наконец, в \cite{AShest12} были описаны тернарные
дифференцирования сепарабельных йордановых алгебр.


В данной работе описываются обобщенные и тернарные
дифференцирования конечномерных полупростых йордановых супералгебр
над алгебраически замкнутым полем характеристики 0, а также для
простых йордановых супералгебр с полупростой четной частью над
алгебраически замкнутым полем произвольной характеристики не 2.
Показано, что за исключением единственного особого случая, любое
тернарное дифференцирование покомпонентно является суммой
обыкновенного дифференцирования и элементов центроида.

\medskip

Пусть $P$ --- поле. Супералгебра $A=A_0+A_1$
--- это  ${\rm Z_2}$-градуированная $P$-алгебра, т.\,е.
$$
A_0^2\subseteq A_0, \  A_1^2\subseteq A_0, \ A_1A_0\subseteq A_1,
\ A_0A_1\subseteq A_1.
$$
Пространство $A_0$ ($A_1$) называется четной (нечетной) частью
супералгебры $A$. Элементы множества $A_0\cup A_1$ называются
однородными. Выражение $\widehat{x}$, где $x\in A_0\cup A_1$,
означает индекс четности однородного элемента $x$:
\,$\widehat{x}=0$, если $x\in A_0$ ($x$
--- четный) и \,$\widehat{x}=1$, если $x\in A_1$ ($x$ --- нечетный).

Линейное отображение $\phi$ супералгебры $A$ называется однородным
(индекса $j$), если $\phi(A_i)\subseteq A_{i+j(mod\ 2)}$,
$i,j=0,1$. Линейное однородное отображение индекса 0(1) называется
четным (соответственно, нечетным).

Множество однородных линейных отображений 
$$
C(A)=\{\,\phi\in{\rm End}(A)\ |
\ \phi(xy)=\phi(x)y=(-1)^{\widehat{x}\widehat{\phi}}x\phi(y)\ \
\forall\,x,y \in A_0\cup A_1\,\}
$$
называется {\it суперцентроидом} $A$.

Однородное линейное отображение $D\in {\rm End}(A)$ называется
{\it супердифференцированием} $A$, если для любых однородных $x,y
\in A$ выполняется равенство
$$
D(xy)=D(x)y+(-1)^{\widehat{x}\widehat{D}}xD(y).
$$
Множество всех супердифференцирований $A$ обозначается как $Der(A)$.


Далее в случаях, когда речь идет о структурах, заданных на
супералгебре, если не оговорено противное, будем подразумевать их
однородными и для краткости будем опускать приставку 'супер',
т.\,е. говоря о различного рода дифференцированиях, центроидах,
подалгебрах, идеалах, и т.\,п. в супералгебре, имеем в виду
супердифференцирования, суперцентроиды, и т.\,д.


Тройка $\Delta=(D,F,G)$ однородных линейных отображений $D,F,G \in
{\rm End}(A) $ называется {\it тернарным
(супер)дифференцированием} cупералгебры $A$, если для любых
однородных $x,y \in A$ выполняется равенство
\begin{eqnarray}\label{komp}
D(xy)=F(x)y+(-1)^{\widehat{x}\widehat{G}}xG(y).
\end{eqnarray}
Как видно из этого определения, тройка вида $(D,D,D)$, где $D\in
Der(A)$ --- дифференцирование $\mbox{супералгебры}\ A$, является
тернарным дифференцированием. Аналогично тройка
$(\phi,\delta\phi,\delta\phi)$ --- тернарное дифференцирование для
$\delta\mbox{-дифференцирования}\ \phi$.

\medskip

\texttt{ЗАМЕЧАНИЕ.} Отметим, что все три отображения $D,F,G$ имеют
одну четность. Действительно, рассмотрим четности суммы и
слагаемых в (1):
$$\widehat{D(xy)}=\widehat{D}+\widehat{xy}=\widehat{D}+\widehat{x}+\widehat{y},$$
$$\widehat{F(x)y}=\widehat{F(x)}+\widehat{y}=\widehat{F}+\widehat{x}+\widehat{y},$$
$$\widehat{xG(y)}=\widehat{x}+\widehat{G(y)}=\widehat{x}+\widehat{G}+\widehat{y}.$$
Поскольку элементы однородные, то все четности
должны быть одинаковы,
$\widehat{D(xy)}=\widehat{F(x)y}=\widehat{xG(y)}$, значит, отсюда
$\widehat{D}=\widehat{F}=\widehat{G}$. Таким образом, можно
говорить и о четности всего тернарного дифференцирования
$\Delta=(D,F,G)$.

\medskip

Множество всех тернарных дифференцирований супералгебры $A$ будем
обозначать $TDer(A)$. Нетрудно убедиться, что оно образует
супералгебру Ли относительно покомпонентных операций векторного
пространства и коммутаторного умножения \nopagebreak
$$ [(D,F,G),(D',F',G')] = ([D,D'],[F,F'],[G,G']), $$
где $[X,Y]=XY-(-1)^{\widehat{X}\,\widehat{Y}}YX$.

Первую компоненту $D$ в тройке тернарного дифференцирования
$(D,F,G)$ будем называть главной, а также дадим ей определение ---
{\em обобщенное дифференцирование}. Множество всех обобщенных
дифференцирований алгебры $A$ будем обозначать $GDer(A)$. Таким
образом,
$$GDer(A)= \{\,D\in{\rm End}(A) \ |\ \exists\,F,G\in{\rm End}(A):\, (D,F,G)\in TDer(A)\,\}, $$
т.\,е. это есть проекция множества тернарных дифференцирований
$TDer(A)$.

Очевидно, что $GDer(A)$, таким же образом, как и $TDer(A)$,
является супералгеброй Ли и образует тем самым подалгебру в $({\rm
End}(A)_0+{\rm End}(A)_1)^{(-)}$, где ${\rm End}(A)_0$  (${\rm
End}(A)_1 )\text{ --- }$ четные (нечетные) линейные отображения
супералгебры $A$. В свою очередь, ясно, что все обыкновенные
дифференцирования являются обобщенными, и их множество $Der(A)$, в
свою очередь, образует подалгебру в $GDer(A)$ относительно тех же
операций
--- обычного векторного пространства и коммутаторного умножения операторов.

\medskip

Заметим, что тройки $\Delta=(D,F,G)$ отображений следующего вида,
очевидно, являются тернарными дифференцированиями в любой
супералгебре:
\begin{eqnarray}\label{komp2}
(D,F,G)=(\varphi + D^0,\,\chi + D^0,\,\psi + D^0),
\end{eqnarray}
где $\varphi,\chi,\psi\in C(A)$ --- произвольные элементы
центроида cупералгебры $A$, удовлетворяющие условию
\,$\varphi=\chi+\psi$, \,а $D^0\in Der(A)$
--- любое обыкновенное дифференцирование в $A$, т.\,е.
$$D^0(xy)=D^0(x)y+(-1)^{\widehat{x}\widehat{D}}xD^0(y),$$
при этом четность у всех одна и та же,
$\widehat{\varphi}=\widehat{\chi}=\widehat{\psi}=\widehat{D^0}=\widehat{\Delta}$.
Назовем все такие тернарные дифференцирования, т.\,е. вида
(\ref{komp2}) --- {\it стандартными}.



Отметим, что для полупервичных супералгебр (за рамки которых в
дальнейшем выходить не будем) нечетная часть центроида равна нулю,
поэтому в нечетном случае ($\widehat{\Delta}=1$) имеем
$\varphi=\chi=\psi=0$, и тогда стандартное тернарное
дифференцирование имеет вид
$$
(D,F,G)=(D^1, D^1, D^1),
$$
где $D^1$ --- обыкновенное нечетное дифференцирование в $A$,
т.\,е.
$$D^1(xy)=D^1(x)y+(-1)^{\widehat{x}}xD^1(y).$$


Соответственно, определим стандартные обобщенные дифференцирования
как главные компоненты стандартных тернарных дифференцирований ---
это будут все преобразования вида
\begin{eqnarray}\label{komp3}
D=\varphi + D^0,
\end{eqnarray}
где \,$\varphi$ --- произвольный элемент центроида, \,а $D^0$
--- любое обыкновенное дифференцирование $\mbox{в }A$, при этом
четность их одна и та же,
$\widehat{\varphi}=\widehat{D^0}=\widehat{D}$.
И также отметим, что в полупервичных супералгебрах в нечетном
случае обобщенного дифференцирования ($\widehat{D}=1$)
имеем
$$
 D=D^1
$$
--- обыкновенное нечетное дифференцирование $\mbox{в }A$.

\medskip

Далее в настоящей работе мы покажем, что для основных классов
простых йордановых супералгебр, за одним исключением, стандартные
дифференцирования, приведенные выше в (\ref{komp2})-(\ref{komp3}),
составляют все множество тернарных/обобщенных дифференцирований.


\medskip

Пусть  $G$ --- алгебра Грассмана над $F$ с единицей $1$, заданная
образующими $1,\xi_{1},\ldots ,\xi_{n},\ldots $ и определяющими
соотношениями: $\xi_{i}^{2}=0,\ \xi_{i}\xi_{j}=-\xi_{j}\xi_{i}.$
Произведения
$$1, \xi_{i_{1}}\xi_{i_{2}}\ldots \xi_{i_{k}}\ \ (i_{1}<i_{2}< \ldots
<i_{k})$$ образуют базис алгебры $G$ над $F$. Обозначим через
$G_{0}$ и $G_{1}$ подпространства, порожденные соответственно
произведениями четной и нечетной длины. Тогда $G = G_{0}+G_{1}$
--- ${\rm Z_2}$-градуированная $P$-алгебра.

Пусть теперь $A=A_{0}+A_{1}$ --- произвольная супералгебра над
$P$. Рассмотрим тензорное произведение $P$-алгебр $G \otimes A$.
Его подалгебра
\begin{eqnarray*}G(A)&=&G_{0} \otimes
A_{0} + G_{1} \otimes A_{1}
\end{eqnarray*}
называется грассмановой оболочкой супералгебры $A$.

Пусть $P$ --- поле характеристики не равной 2.  Алгебра
$J=J_{0}\oplus J_{1}$ называется йордановой супералгеброй, если ее
грассманова оболочка $G(J)$ является йордановой алгеброй, т.\,е. в
$G(J)$ справедливы тождества:
$$xy=yx,\ \ (x^2y)x=x^2(yx).$$


\medskip

Основным результатом о строении простых конечномерных йордановых
супералгебр является следующие классификационные теоремы (см.
также \cite{IShest97}):

\medskip


\textbf{Теорема\cite{Kantor89},\cite{Kac77}.} Всякая конечномерная
простая нетривиальная (т.\,е. с ненулевой нечетной частью)
йорданова супералгебра $J$ над алгебраически замкнутым полем $P$
характеристики 0 изоморфна одной из следующих супералгебр:

i) $M_{m,n}^{(+)}(P),\ m,n>0$;

ii) $osp(n,2m),\ m,n>0$;

iii) $P(n)\ (\,= trp(n,n)$, см.\cite{IShest97}\,$),\ n>1$;

iv) $Q(n)\ (\,= M_n[u]^{(+)}$, см.\cite{IShest97}\,$),\ n>1$;

v) супералгебра билинейной формы $J(V,f)$;

vi) супералгебра Капланского $K_3$;

vii) $D_t$, $t\in K,\ t\neq 0$;

viii) супералгебра Каца $K_{10}$;

ix) супералгебра грассмановых скобок Пуассона $J(\Gamma_n)$.

\medskip

\textbf{Теорема\cite{RacZelm}.} Всякая конечномерная простая
нетривиальная (т.\,е. с ненулевой нечетной частью) йорданова
супералгебра $J$ c полупростой четной частью над алгебраически
замкнутым полем $P$ характеристики не равной 2 изоморфна одной из
следующих супералгебр:

i)-viii) из предыдущей теоремы;

x) вырожденная супералгебра Каца $K_9,\ char\,P=3$;

xi)  $H_3(B(1,2))$ (см.\cite{IShest97}), $char\,P=3$;

xii) $H_3(B(4,2))$ (см.\cite{IShest97}), $char\,P=3$.



\medskip

Отметим, что из приведенных супералгебр во всех случаях, кроме
(vi) --- это алгебры с единицей, а во всех случаях, кроме (ix)
--- алгебры с полупростой четной частью, поэтому далее разделим
изложение на части --- сначала рассмотрим задачу для любых
унитальных супералгебр с полупростой четной частью, затем отдельно
супералгебру Капланского $K_3$ и отдельно супералгебру
грассмановых скобок $J(\Gamma_n)$. Еще, выделим случай (v') ---
\begin{eqnarray}\label{komp v'}
J(V,f)=Pe+V,\ \text{где}\ V=V_0+V_1\ \,\text{при}\ V_0=0,\, V_1=V,
\end{eqnarray}
--- супералгебру билинейной формы с одномерной четной частью,
т.к. забегая вперед, скажем, что для этого случая результат будет
отличаться от основного.


\medskip

{\bf Лемма 1.} {\em Пусть $B$ --- произвольная коммутативная
унитальная супералгебра. Тогда линейное преобразование $D$ будет
обобщенным дифференцированием $B$ тогда и только тогда, когда в
ней $D$ удовлетворяет равенству: }
\begin{eqnarray}\label{komp4}
D(xy) =
D(x)y+(-1)^{\widehat{x}\widehat{D}}xD(y)-(-1)^{\widehat{x}\widehat{D}}((xc)y+x(cy))\
\ \ \forall \,x,y \in B,
\end{eqnarray}
{\em где $c=D(1)/2$.}

\medskip

{\bf Доказательство} Пусть $D\in GDer(B)$ --- обобщенное и
$(D,F,G) \in TDer(B)$
--- связанное с $D$ тернарное дифференцирование, т.\,е.
выполняется определяющее равенство (\ref{komp})
$$D(xy)=F(x)y+(-1)^{\widehat{x}\widehat{G}}xG(y)=F(x)y+(-1)^{\widehat{x}\widehat{D}}xG(y).$$

Положим в равенстве (1) $y=1$. Тогда
$$ D(x)=F(x)+
(-1)^{\widehat{x}\widehat{D}}xG(1).$$ Обозначив $G(1)=g$ (отметим
при этом, что тогда $\widehat{g}=\widehat{D}$), имеем
$$F(x)=D(x)-(-1)^{\widehat{x}\widehat{D}}xg=D(x)-(-1)^{\widehat{x}\widehat{D}}(-1)^{\widehat{x}\widehat{D}}gx=D-gx,$$
или $F=D-L_g$, \,где $L_a$ --- оператор левого умножения на элемент
$a\in J$, т.\,е. $L_a(x)=ax$. C другой стороны,
$$x=1 \ \ \Longrightarrow \ \ D(y)=F(1)y + G(y),$$
откуда, обозначив $F(1)=f$ (при этом $\widehat{f}=\widehat{D}$),
имеем $G(y)=D(y)-fy$, \,или \,$G=D-L_f$.

Подставляя  выражения $F,G$ снова в (\ref{komp}), получаем
$$
D(xy) = D(x)y + (-1)^{\widehat{x}\widehat{D}}xD(y) -(gx)y -
(-1)^{\widehat{x}\widehat{D}}x(fy) =
$$
\begin{eqnarray}\label{komp5}
= D(x)y + (-1)^{\widehat{x}\widehat{D}}xD(y)
-(-1)^{\widehat{x}\widehat{D}}(xg)y -
(-1)^{\widehat{x}\widehat{D}}x(fy)
\end{eqnarray}
С другой стороны,
$$ D(yx)=D(y)x + (-1)^{\widehat{y}\widehat{D}}yD(x) -(gy)x -(-1)^{\widehat{y}\widehat{D}}y(fx)= $$
$$ =(-1)^{\widehat{x}\widehat{D(y})}xD(y)+(-1)^{\widehat{y}\widehat{D}}(-1)^{\widehat{y}\widehat{D(x)}}D(x)y
-(-1)^{\widehat{x}\,\widehat{gy}}x(gy)-(-1)^{\widehat{y}\widehat{D}}(-1)^{\widehat{y}\widehat{fx}}(fx)y=$$

$$=(-1)^{\widehat{x}(\widehat{y}+\widehat{D})}xD(y) +
(-1)^{\widehat{y}\widehat{D}}(-1)^{\widehat{y}(\widehat{x}+\widehat{D})}D(x)y
- (-1)^{\widehat{x}(\widehat{y}+\widehat{D})}x(gy) -
(-1)^{\widehat{y}\widehat{D}}(-1)^{\widehat{y}(\widehat{x}+\widehat{D})}(-1)^{\widehat{x}\widehat{f}}(xf)y=
$$
$$
=(-1)^{\widehat{x}\widehat{y}+\widehat{x}\widehat{D}}xD(y) +
(-1)^{\widehat{x}\widehat{y}}D(x)y -
(-1)^{\widehat{x}\widehat{y}+\widehat{x}\widehat{D}}x(gy) -
(-1)^{\widehat{x}\widehat{y}+\widehat{x}\widehat{D}}(xf)y ;
$$
В то же время, в силу суперкоммутативности умножения
$$
D(yx)=(-1)^{\widehat{x}\widehat{y}}D(xy)
=(-1)^{\widehat{x}\widehat{y}}D(x)y +
(-1)^{\widehat{x}\widehat{y}+\widehat{x}\widehat{D}}xD(y) -
(-1)^{\widehat{x}\widehat{y}+\widehat{x}\widehat{D}}(xg)y -
(-1)^{\widehat{x}\widehat{y}+\widehat{x}\widehat{D}}x(fy),
$$
откуда
$$(xg)y + x(fy) = x(gy) + (xf)y, $$
$$ x(fy-gy)=(xf-xg)y, $$
$$ x((f-g)y)=(x(f-g))y, $$
или, обозначив $f-g=w$, запишем в форме ассоциатора:
\begin{eqnarray}\label{komp6}
(x,w,y)=0
\end{eqnarray}
--- это равенство выполняется для любых $x,y \in B$.

\medskip

Имеем:
$$\left\{\aligned f+g=D(1) \\f-g=w \endaligned \right.\ , $$
откуда $f=c+w/2$, $g=c-w/2$, \ где\, $c=D(1)/2$.

\medskip

Теперь, возвращаясь к (\ref{komp5}), запишем:
$$D(xy)=D(x)y + (-1)^{\widehat{x}\widehat{D}}xD(y) -
(-1)^{\widehat{x}\widehat{D}}[(x(c-w/2))y - x((c+w/2)y)]=
$$
$$
= D(x)y + (-1)^{\widehat{x}\widehat{D}}xD(y)
-(-1)^{\widehat{x}\widehat{D}}[(xc)y -x(cy) + \frac{1}{2}((xw)y
-x(wy))],
$$
откуда в силу (\ref{komp6}) получаем равенство (\ref{komp4}).

Нетрудно проверить, что верно и обратное. Еcли $D$ ---
удовлетворяет равенству (\ref{komp4}), тогда взяв любые $f,g \in
B$ такие, что $f+g=D(1)$, $(B,f-g,B)=0$, (например,
$f=g=c=D(1)/2$), и положив $F=D-L_g$, $G=D-L_f$, получим тернарное
дифференцирование $(D,F,G)$.
Лемма доказана.


\medskip

\texttt{ЗAМЕЧАНИЕ 1.} Из леммы 1 следует, что обобщенное
дифференцирование $D$ в унитальной коммутативной
$\mbox{супералгебре}\ B$ \,является обыкновенным
дифференцированием в том и только в том случае, когда $D(1)=2c=0$,
т.\,е. как и в обычных алгебрах

$$Der(B)=\{\,D\in GDer (B)\ |\ D(1)=0\,\}.$$

\texttt{ЗAМЕЧАНИЕ 2.} Из предыдущего также следует, что обобщенное
дифференцирование $D$ в унитальной коммутативной
$\mbox{супералгебре}\ B$ \,является стандартным в том и только в
том случае, когда $D(1)=\alpha\cdot 1\in Z(B)$ --- элемент из
центра $B$.


\medskip

Помимо прочего, в лемме показано, что всякое тернарное
дифференцирование в унитальной коммутативной cупералгебре $B$
имеет вид \ $\Delta=(D,\,D-L_g,\,D-L_f)$, $D\in GDer\,(B)$ ---
таким образом все тернарные дифференцирования сводятся к
обобщенным.

\medskip

{\bf Теорема 1.} {\em Пусть $J$ --- конечномерная простая
унитальная йорданова супералгебра с полупростой четной частью
над алгебраически замкнутым полем $P$, отличная от супералгебры
(\ref{komp v'}). Тогда все oбобщенные дифференцирования в $J$
являются стандартными. т.\,е. вида (\ref{komp3})}.


\medskip

Для доказательства теоремы приведем еще некоторые нужные нам
понятия и свойства.

\medskip

Пусть $J$ --- йорданова (супер)алгебра с единицей $1=\sum_{i=1}^n
e_i$, разлагающейся в сумму попарно ортогональных идемпотентов
$e_i^2=e_i\in J$. Тогда (см.\cite{Jacob68}) имеет место разложение
в прямую сумму подпространств $J=\bigoplus_{(1\leq i\leq j\leq n)}
J_{ij}$, \,где

$$ J_{ii}=\{ x\in J\ |\ xe_i=x \}, \ \ J_{ij}=\{ x\in J\ |\ xe_i=xe_j=\frac{1}{2}x\},$$
причем подпространства $J_{ij}$ связаны соотношениями

$$ J_{ii}^2\subseteq J_{ii}\,, \ \ J_{ij}J_{ii}\subseteq J_{ij}\,, \ \ J_{ij}^2\subseteq J_{ii}+J_{jj}\,,$$

$$ J_{ii}J_{jj}=J_{ii}J_{jk}=J_{ij}J_{kl}=(0), \ \ J_{ij}J_{jk}\subseteq J_{ik}\,,$$ где индексы $i,j,k,l$
все различны. Приведенное разложение называется {\em пирсовским
разложением} $J$ относительно системы идемпотентов
$e_1,\ldots,e_n$.


\medskip

{\bf Доказательство.} Сразу зафиксируем, что в условиях данной
теоремы имеет место все вышеперечисленное, причем отметим, что
систему ортогональных идемпотентов в $J$ --- $e_1,\ldots, e_n$,
$\sum_{i=1}^n e_i = 1$ --- можно подобрать так, что все
одноименные подпространства одномерны, т.\,е. $J_{ii}=P\cdot e_i$
для всех $i$. Кратко опишем системы идемпотентов и соответствующие
им подпространства в каждом случае --- здесь это будут все
супералгебры, кроме (vi),(ix) и выделенного случая (v').

i) $J=M_{m,n}^{(+)}$; здесь идемпотенты $e_i=e_{ii}$ --- матричные
единицы, и нетрудно проверить, что соответствующее $e_i$
пирсовское подпространство $J_{ii}=Pe_{ii}$, \ $i=1,\ldots,m+n$.

ii) $J=osp(n,2m)$; \ $e_i=e_{ii},\ \text{при}\ i\leq n$, \
$e_i=e_{ii}+e_{m+i,m+i},\ \text{при}\ i>n$, \ $i=1,\ldots,n$.

iii) $J=trp(n,n)$; \ $e_i=e_{ii}+e_{n+i,n+i}$, \
$J_{ii}=P(e_{ii}+e_{n+i,n+i})$, \ $i=1,\ldots,n$.

iv) $J=M_n[u]^{(+)}$; \ $e_i=e_{ii}+e_{n+i,n+i}$,
$J_{ii}=P(e_{ii}+e_{n+i,n+i})$, \ $i=1,\ldots,n$.

v) $J=J(V,f)=Pe+V$,\, $V=V_0+V_1$, причем $V_0\neq 0$. Тогда
$\exists\,v\in V_0$: $v^2=f(v,v)=1$. Взяв $e_1=(1+v)/2$,
$e_2=(1-v)/2$, получим $e_1+e_2=1$, $e_1e_2=0$, $e_i^2=1$; при
этом $J_{11}=Pe_1$, $J_{22}=Pe_2$, $J_{12}=(Pv)^{\bot}$,
$J=J_{11}+J_{12}+J_{22}$.

vii) $J=D_t=(D_t)_0+(D_t)_1$, где $(D_t)_0=Pe_1+Pe_2$, \,
$(D_t)_1=Px+Py$, где $e_i^2=1$, $e_1e_2=0$, $e_i x=x/2$, $e_i
y=y/2$; здесь базисные элементы $e_1,e_2$ являются системой
идемпотентов по определению, т.\,к. можно заметить, что
$e_1+e_2=1$ в $J$; и также $J_{11}=Pe_1$, $J_{22}=Pe_2$.

viii) $J=K_{10}$, \ $J_0=A_1\oplus A_2$, \ $J_1=Pu+Pv+Pz+Pw$,
\,где (см.\cite{IShest97}) $A_1=Pe_1+Puz+Puw+Pvz+Pvw$, \
$A_2=Pe_2$, \ $e_i^2=e_i$, \ $e_1x=x\ \ \forall x\in A_1$, \
$e_ix= x/2\ \ \forall x\in J_1$. Нетрудно видеть, что $e_1+e_2=1$,
и  $e_1,e_2$
--- система идемпотентов, однако подпространство, сооответствующее $e_1$, не
одномерно. Представим в виде $J_0=Pe_2+J(V,f),\ где
V=L(uz,uw,vz,vw)$. Тогда, как и \text{в случае (v)},
$e_1=e_1'+e_2''$, здесь уже соответсвующие подространтсва
одномерны. Имеем $1=e_1'+e_1''+e_2$, таким образом положив
$\tilde{e_1}=e_1',\ \tilde{e_2}=e_1''\ \tilde{e_3}=e_2$, получим
систему идемпотентов с одномерными соответствующими
подпространствами $J_{ii}=P\tilde{e_i}$.

x) $J=K_{9}$, --- аналогично.

xi)-xii) $J=H_3(B)$, где $B\in\{B(1,2), B(4,2)\}$; как и для
матричных супералгебр, системой идемпотентов являются матричные
единицы $e_i=e_{ii}$, и соответствено, $J_{ii}=Pe_i$, $i=1,2,3$.

\medskip

Пусть теперь $D\in GDer(J)$
--- обобщенное, и $(D,F,G) \in TDer(J)$
--- связанное с $D$ тернарное дифференцирование, т.\,е.
выполняется определяющее равенство (\ref{komp})
$$D(xy)=F(x)y+(-1)^{\widehat{x}\widehat{G}}xG(y)$$


Случай 1. $D$ --- четное, т.\,е. $\widehat{D}=0$. Тогда основное
равенство (1) \mbox {для $D$:}
$$ D(xy)=F(x)y + xG(y),$$
т.\,е. тернарное и соответственно, обобщенное дифференцирование в
этом случае действует, как в обычных алгебрах. Рассмотрим сужение
$D$ на четную часть $J_0$: $D|J_0=D_0\in GDer(J_0)$ --- обобщенное
дифференцирование обычной йордановой алгебры $J_0$. В силу
результата \cite{AShest12} для обычных йордановых алгебр, $D_0$
является стандартным (в \cite{AShest12} именовалось
"тривиальным"), т.\,е. $D_0=\alpha I+D^0_0$, где $D^0_0\in
Der(J_0)$
--- обычное дифференцирование обычной алгебры $J_0$, a
\,$\alpha\in Z(J_0)$.

\medskip

а) Если алгебра $J_0$ --- проста, то $\alpha\in P$ --- скаляр.
Тогда уже на всей супералгебре $J$ рассмотрим $D-\alpha I=D^0\in
GDer(J)$
--- тоже обобщенное диференцирование, при этом
$D^0(1)=D(1)-\alpha\cdot 1=D_0(1)-\alpha=0$, значит, согласно
замечанию после леммы 1, $D^0\in Der(J)$, и таким образом,
$D=\alpha I + D^0$ --- стандартное четное обобщенное
дифференцирование.

\medskip

б) Если алгебра $J_0$ --- непроста, то (\cite{RacZelm})
$J_0=A_1\oplus A_2$, где $A_i$ --- простые алгебры. Тогда
$D|A_i=D_i\in GDer(A_i)$
--- обобщенные дифференцирования в $A_i$ и они стандартны,
что согласно \text{замечанию 2} после \text{леммы 1} значит
$D_i(e_i)=\alpha_ie_i$, где $e_i$
--- единицы в $A_i$, $\alpha_i\in P$ --- скаляры.
Покажем, что $\alpha_1=\alpha_2=\alpha$, после чего
можно будет действовать как в \text{п. а)}.
Заметим, что $e_1+e_2=1$ --- система ортогональных идемпотентов
$\text{в } J_0$ и также $\text{в }J$,
--- таким образом, имеет место пирсовское разложение
$J=J_{11}+J_{12}+J_{22}$, причем (см.\cite{RacZelm}) смешанная
компонента $J_{12}=J_1$ --- совпадает с нечетной частью исходной
$\text{супералгебры }J$. Имеем
$$
D(e_i)=\alpha_ie_i, \ \,e_i^2=e_i, \ \,e_1e_2=0, \ \,e_1+e_2=1, \
\,e_1x=e_2x=\frac{1}{2}x\, \ \forall\,x\in J_1.
$$

Применяя $D$ к произвольному $h\in J_1$ (при этом, т.\,к.
$\widehat{D}=0$, то и $D(h)\in J_1=J_{12}$), получаем
$$
D(h) = 2D(e_1h) = 2D(e_1)h + 2e_1D(h) - (e_1
 D(1))h  - e_1 (D(1)h) =
$$
$$ = 2\alpha_1 e_1h + D(h) - (e_1(\alpha_1
e_1 + \alpha_2 e_2))h  -  e_1((\alpha_1 e_1 + \alpha_2 e_2) h) =
$$
$$
= \alpha_1 h + D(h) - \alpha_1 e_1h - \frac{1}{2}e_1(\alpha_1 h +
\alpha_2 h)\, =\,D(h) + \frac{1}{4}(\alpha_1 - \alpha_2)h ,
$$
откуда $\alpha_1-\alpha_2 = 0$, $\alpha_1 =\alpha_2 =\alpha\in P$.
Тогда $D(1)=\alpha\cdot 1$ --- свели к п. а).



\medskip

Случай 2. $D$ --- нечетное, т.\,е. $\widehat{D}=1$. Тогда основное
равенство (1) \mbox {для $D$:}
$$ D(xy)=F(x)y + (-1)^{\widehat{x}}xG(y)$$

Пусть $1=\sum_{i=1}^ne_i$, где $e_i$ --- ортогональные идемпотенты
из $J_0$. Рассмотрим 
пирсовское
 разложение $J$ относительно $e_1,\ldots, e_n$, тогда, применив \mbox{лемму 1},
получаем:
$$
D(e_i)=2D(e_i)e_i-2(e_i c)e_i \in J\cdot e_i = J_{ii} + \sum_{j\ne
i} J_{ij},
$$
где $J_{ii},J_{ij}$ --- подпространства пирсовского разложения.
Таким образом,
$$
D(e_i) = \alpha_i  e_i + \sum_{j\ne i} d_{ij}, \ \ \alpha_i\in P,
\ d_{ij}\in J_{ij},
$$
и соответственно,
$$
D(1) = \sum_{i=1}^n D(e_i)=\sum_{i=1}^n \alpha_i  e_i +
\sum_{i=1}^n\sum_{j\ne i} d_{ij}=\sum_{i=1}^n \alpha_i e_i +
\sum_{i\ne j} d_{ij}=2c.
$$

С другой стороны, $ D(e_i)=2D(e_i)e_i-2(e_i c)e_i = $
$$
=2\alpha_i e_i+\sum_{j\ne i} d_{ij} - (e_i(\sum_{k=1}^n\alpha_k
e_k + \sum_{k\ne j}d_{kj}))e_i =
$$
$$
= 2\alpha_i e_i+\sum_{j\ne i} d_{ij} - (\alpha_i e_i +
\frac{1}{2}\sum_{j\ne i}(d_{ij}+d_{ji}))e_i = \alpha_i e_i +
\sum_{j\ne i} d_{ij} - \frac{1}{4}\sum_{j\ne i}(d_{ij}+d_{ji}),
$$
значит, отсюда следует, что $d_{ij}+d_{ji}=0$ для всех $i\ne j$.

Тогда
$$2c=D(1)=\sum_{i=1}^n \alpha_i e_i + \sum_{i\ne j} d_{ij} =
\sum_{i=1}^n \alpha_i e_i \in J_{0}.$$

В то же время, т.\,к. $\widehat{D}=1$, то $D(1)\in J_1$.

Значит, остается $D(1)=0=2c$. Тогда равенство $(\ref{komp4})$ из
леммы 1 для $D$ будет выглядеть:
$$D(xy)=D(x)y+(-1)^{\widehat{x}}xD(y),$$
и таким образом, $D$ --- обыкновенное нечетное дифференцирование.
Теорема доказана.



\medskip

В работе о тернарных дифференированиях обычных йордановых алгебр
для описания двух оставшихся компонент $F$ и $G$ тернарного
дифференцирования мы использовали результат McCrimmon-Ng
\cite{McCrimmon82} о среднем центре полупервичных йордановых
алгебр. Введем несколько определений и понятий, чтобы привести
изложение в соответствие с терминологией McCrimmon-Ng.

\medskip

Идеал $I$ в йордановой алгебре называется {\em тривиальным}, если
$U_I I=0$, где оператор $U_x=2L_x^2-L_{x^2}$.

\medskip
Йорданова алгебра называется \emph{полупервичной}, если она не
содержит тривиальных идеалов; и кроме того, называется
\emph{$D$-полупервичной}, если она не содержит триваильных
идеалов, замкнутых относительно всех дифференцирований алгебры.

\medskip

{\bf Теорема (McCrimmon-Ng \cite{McCrimmon82}).} {\em Для всякой
$D$-полупервичной йордановой алгебры $J$ средний ассоциативный
центр алгебры $J$, т.\,е. множество элементов
$$ W(J) = \{\,w\in J\ |\ (x,w,y)=0\ \ \forall x,y\in J \,\} $$
совпадает с центром алгебры $J$, $W(J)= Z(J)$.}

\medskip

Нам нужно теперь обобщить этот результат для супералгебр. Для
этого докажем еще несколько вспомогательных утверждений.

\medskip

{\bf Лемма 2.} \emph{ 1) Тривиальный в определении McCrimmon-Ng
идеал --- означает в точности такой идеал $I$ йордановой алгебры
$J$, что $I^3=0$. 2) Условие полупервичности в смысле McCrimmon-Ng
для йордановой алгебры равносильно условию отсутствия в ней
ненулевых нильпотентных идеалов.}

\medskip

{\bf Доказательство.} 1) $I^3=0\ \Longrightarrow\ U_I(I)=0$ ---
просто по определению. Покажем обратно. Пусть $I$ --- тривиальный
по McCrimmon-Ng идеал, т.\,е. $U_I(I)=0$. Возьмем произвольные
элементы $a,b,c\in I$. Рассмотрим
$U_{a,b}(c)=U_{a+b}(c)-U_a(c)-U_b(c)=0$. Раскрывая и приводя
подобные, а также в силу коммутативности умножения, получаем в
итоге $U_{a,b}(c)=(ac)b+a(bc)-(ab)c=0$. Аналогично рассмоотрим
$U_{a,c}(b)=-(ac)b+a(bc)+(ab)c=0$,
$U_{b,c}(a)=(ac)b-a(bc)+(ab)c=0$. Обозначив $(ac)b=x$, $a(bc)=y$,
$(ab)c=z$, получаем однородную систему линейных уравнений
относительно $x,y,z$, определитель которой оказывается не равен 0,
--- значит, она имеет только тривиальное решение соответственно
$(ac)b=a(bc)=(ab)c=0$. Поскольку элементы $a,b,c\in I$ брались
произвольные, то значит, $I^3=0$.

2) При отсутствии нильпотентных идеалов условие полупервичности
McCrimmon-Ng в связи с доказанным п.1) выполнено по определению.
Обратно --- предположим есть нильпотентный идеал $I^n=0$; пусть
$3^k\leq n \leq 3^{k+1}$. Тогда $I'=(..((I\underbrace{^3)^3)\ldots
\ldots )^3}_k$
--- тоже является идеалом, при этом $(I')^3=0$. Лемма доказана.

\medskip

В связи с доказанным, и ввиду отсутствия определения полупервичности
у McCrimmon-Ng для супералгебр, в дальнейшем будет иметь смысл
ввести обобщенное определение для любых алгебр и супералгебр.

\medskip

Алгебра (супералгебра) называется \emph{сильно полупервичной},
если она не содержит нильпотентных идеалов; и соответственно,
\emph{$D$-сильно полупервичной}, если она не содержит
нильпотентных идеалов, замкнутых относительно всех
дифференцирований алгебры.

\medskip

\texttt{ЗAМЕЧАНИЕ.} Для обычных йордановых алгебр введенное выше
понятие сильной полупервичности и $D$-полупервичности совпадает с
определением по McCrimmon-Ng, что следует из леммы 2.


\medskip

{\bf Предложение.} {\em Если $J$ --- сильно полупервичная
супералгебра, тогда ее грассманова оболочка $G(J)$ будет
$D$-сильно полупервична.}

\medskip

{\bf Доказательство.} Пусть $I\lhd G(J)$ --- нильпотентный идеал,
замкнутый относительно всех дифференцирований, т.\,е. $I^n=0$ \ и
\ $D(I)\subseteq I$ \ \,$\forall\,D\in Der G(J)$. Покажем что
тогда $I=0$.

\medskip

Предположим, cуществует элемент $0\ne i\in I\subset G(J)$, \,т.\,е.
\,$i=i_0+i_1$,
$$
i_0=\sum_{k}\xi_k\otimes a_k, \ \ i_1=\sum_{l}\eta_l\otimes b_l, \
\ \text{где}\ a_k\in J_0,\ \xi_k\in G_0,\ b_l\in J_1,\ \eta_l\in
G_1.
$$
Перенумеруем индексы, чтобы $\xi_1$ --- одночлен минимальной
длины среди $\{\xi_k\}$, $\eta_1$
--- соответственно, среди $\{\eta_l\}$. Пусть $\zeta=e_{j_1}\cdots e_{j_p}$ ---
минимальный одночлен  из \,${\xi_1,\eta_1}$, т.\,е. таким образом
минимальный из всех одночленов $\{\xi_k\}$,$\{\eta_l\}$. Тогда в
этом случае можно построить одночлен $h=e_{j_{p+1}}\cdots e_{j_q}$
--- содержащий все переменные, входящие в $\{\xi_k\}$,$\{\eta_l\}$
кроме $\zeta$, \,т.\,е. $e_{j_r}\notin\zeta,\ r=p+1,\ldots,q$.
Сделаем $h$ четным, домножив или не домножив, в зависимости от
надобности, на еще один базисный элемент $e_{j_{q+1}}$.

В силу того, что $I\lhd G(J)$, имеем $I\ni i\cdot(h\otimes 1)=$
$$
=\sum_{k}\xi_k h\otimes a_k + \sum_{l}\eta_l h\otimes b_l = \xi_1
h\otimes a_1 + \eta_1 h\otimes b_1 = \zeta h\otimes c,
$$
где $c=a_1$ в случае, если $\zeta=\xi_1$, и $c=b_1$ в случае
$\zeta=\eta_1$.

Таким образом, в результате получили некоторый элемент
$e_{j_1}\cdots e_{j_s}\otimes c\in I$ или, перенумеровав базисные
элементы, будем считать $e_1\cdots e_s\otimes c\in I$.

Рассмотрим в $J$  идеал $I'=id<c>$, порожденный элементом $c$. В
силу условия cильной полупервичности $J$ имеем $(I')^n\ne 0$.
Поэтому $\exists\,c^{(1)},\ldots, c^{(n)}\in I': \ c^{(1)}\cdots
c^{(n)}\ne 0$. \,Имеем

$$
id<c>\ \,\ni\
c^{(k)}=\sum_{p_k=1}^{l_k}(((cx_{p_k1}^{(k)})x_{p_k2}^{(k)})\cdots)x_{p_kq_{p_k}}^{(k)}=
\sum_{p_k=1}^{l_k} u_{p_k}^{(k)}.
$$
Таким образом, $\sum_{p_1=1}^{l_1}
u_{p_1}^{(1)}\cdots\sum_{p_n=1}^{l_n} u_{p_n}^{(n)} \ne 0$,
следовательно, $\exists\,p_1,\ldots p_n: \ u_{p_1}^{(1)}\cdots
u_{p_n}^{(n)}\ne 0$, и значит, без ограничения общности можно
считать, что
$$c^{(k)}=(((cx_1^{(k)})x_2^{(k)})\cdots)x_{m_k}^{(k)}, \ \ k=1,\ldots,n.$$
Наберем для $\{x_i^{(k)}\}$ одночлены $\{\xi_i^{(k)}\in G\}$ с
соответствующей четностью $\widehat{\xi_i^{(k)}}=
\widehat{x_i^{(k)}}$, тогда
$$
I\ni(((e_1\cdots e_s\otimes c) (\xi_1^{(k)}\otimes
x_1^{(k)})\cdots)(\xi_{m_k}^{(k)}\otimes x_{m_k}^{(k)})=
$$
$$
= e_1\cdots e_s \xi_1^{(k)}\cdots
\xi_{m_k}^{(k)}\otimes(((cx_1^{(k)})x_2^{(k)})\cdots)x_{m_k}^{(k)} =
e_1\cdots e_s \xi_1^{(k)}\cdots \xi_{m_k}^{(k)}\otimes c^{(k)} =
c'^{(k)}\in I ;
$$
$\xi_1^{(k)},\ldots\xi_{m_k}^{(k)}$ очевидно, можно подобрать так
чтобы все базисные элементы, содержащиеся в них, продолжали ряд
$e_1,\ldots,e_s$, и тогда имеем $c'^{(k)}=e_1\cdots e_s
e_{s+1}\cdots e_{m_k'}\otimes c^{(k)} \in I$.

Таким образом, для каждого $c^{(k)}$ $\exists\,m_k': \
c'^{(k)}=e_1\cdots e_{m_k'}\otimes c^{(k)}\in I$.

Рассмотрим следующий оператор дифференцирования $D\in Der\,G$,
действующий по правилу $D(e_i)=e_{i+1}$. Cоответственно, тогда
$D(e_{i_1}\cdots e_{i_k})=\sum_{r=1}^k e_{i_1}\cdots
D(e_{i_r})\cdots e_{i_k}$. Нетрудно проверить следующее
примечательное для этого оператора соотношение:
$$
\forall\,k \ \ D^k(e_1e_2\cdots e_k)e_{k+2}\cdots e_{2k+1} =
e_2e_3\cdots e_{2k+1}.
$$
Действительно, по индукции, для $k=1$ --- $D(e_1)e_3=e_2e_3$ ---
верно. Пусть верно для $k-1$, т.\,е. $D^{k-1}(e_1e_2\cdots
e_{k-1})e_{k+1}\cdots e_{2k-1} = e_2e_3\cdots e_{2k-1}$. Тогда для
$k$
$$
D^k(e_1e_2\cdots e_k)e_{k+2}\cdots e_{2k+1} = D^{k-1}(D(e_1\cdots
e_k))e_{k+2}\cdots e_{2k+1}=
$$
$$
=D^{k-1}(e_1\cdots e_{k-1}e_{k+1})e_{k+2}\cdots e_{2k+1}
=\sum_{i+j=k-1}C^i_{k-1} D^i(e_1\cdots
e_{k-1})D^j(e_{k+1})e_{k+2}\cdots e_{2k+1}=
$$
$$
=D^{k-1}(e_1\cdots e_{k-1})e_{k+1}e_{k+2}\cdots e_{2k+1}
=e_2e_3\cdots e_{2k-1}e_{2k}e_{2k+1}
$$
--- получили истинность формулы для $k$.

Кроме того, если данный оператор расширить на $G(J)$ следующим
естественным образом: $\tilde{D}=D\otimes id$, где $id$ ---
тождественное отображение, то можно также убедиться, что
$\tilde{D}(\xi\otimes x) = D(\xi)\otimes x$ и проверить, что
$\tilde{D}\in Der\,G(J)$.

Применим $\tilde{D}^{m_2'}$ к $c'^{(2)}$, тогда
$$
D^{m_2'}(c'^{(2)})(e_{m_2'+2}\cdots e_{3m_2'+1}\otimes
1)=e_2e_3\cdots e_{3m_2'+1}\otimes c^{(2)}\in I.
$$
Продолжая
нужное число раз, можно дойти до $c''^{(2)}=e_{m_1'+1}\cdots
e_{m_2''}\otimes c^{(2)}\in I$ в силу условия инвариантности
относительно дифференцирований.

Аналогично применяя дифференцирование $\tilde{D}$ к $c'^{(k+1)}$
нужное число раз, получим $c''^{(k+1)}=e_{m_k''+1}\cdots
e_{m_{k+1}''}\otimes c^{(k)}\in I$, где по индукции
$c''^{(k)}=e_{m_{k-1}''+1}\cdots e_{m_k''}\otimes c^{(k)}\in I$
($k=2,\ldots, n$).

Таким образом, теперь имеем:
$$
I^n\ni c'^{(1)}c''^{(2)}\cdots c''^{(n)}=e_1\cdots
e_{m_1'}e_{m_1'+1}\cdots e_{m_2''}e_{m_2''+1}\cdots e_{m_k''}\otimes
c^{(1)}c^{(2)}\cdots c^{(n)}\ne 0,
$$
т.\,е. $I^n\ne 0$, что --- противоречие с начальным условием.
Значит, предположение $0\ne i\in I$ --- неверно, и $I=0$. Лемма
доказана.

\medskip

Теперь, собственно, обобщенный результат McCrimmon-Ng:

\medskip

{\bf Лемма 3.} {\em Пусть $J$
---  сильно полупервичная йорданова супералгебра.
Тогда среднее ядро супералгебры $J$, т.\,е. множество элементов
$$ W(J) = \{\,w\in J\ |\ (J,w,J)=0\,\} $$
лежит полностью в центре супералгебры $J$, т.\,е. $W\subseteq
Z(J)$.}

\medskip

{\bf Доказательство.} Пусть $w\in W(J)$ --- произвольный элемент
из среднего ядра cупералгебры $J$, т.\,е.
$$ (x,w,y)=0 \ \ \forall x,y\in J $$
Можно считать, что $w$ --- однородный. Возьмем для $w$ любой
$\xi\in G$ из алгебры Грассмана соответствующей чётности
$\widehat{\xi}=\widehat{w}$, тогда
элемент \ $\xi\otimes w\in G(J)$ будет лежать в центре $G(J)$.
Действительно, для любых произвольных $x^0,y^0\in G(J)$, т.\,е.
$x^0=\xi_1\otimes x$, \,$y^0=\xi_2\otimes y$, \,где \,$x,y\in J$,
\, $\widehat{\xi_1}=\widehat{x}$, $\widehat{\xi_2}=\widehat{y}$,
имеем:

$$ (\xi_1\otimes x,\,\xi\otimes w,\, \xi_2\otimes y)=\xi_1\xi\xi_2\otimes(x,w,y)=0, $$
т.\,е. $(x^0,\,\xi\otimes w,\,y^0)=0 \ \ \forall\,x^0,y^0\in G(J)$
и таким образом, $\xi\otimes w\in W(G(J))$ --- лежит в среднем
ассоциативном цетре  грассмановой оболочки $G(J)$, являющейся
согласно предложению уже $D$-сильно полупервичной йордановой
алгеброй, значит в силу результата \cite{McCrimmon82} McCrimmon-Ng
имеем $\xi\otimes w\in Z(G(J))$
--- является элементом центра для $G(J)$. Причем, поскольку
одночлен $\xi\in G$ для $w$ брали произвольно, то $\xi\otimes w\in
Z(G(J))$ лежит в центре при любом выборе $\xi\in G$.

Сделаем теперь обратный переход и покажем отсюда, что элемент $w$
лежит в центре исходной супералгебры $J$ --- это означает, что
для любых однородных  элементов $u,v\in J$ исходной супералгебры
должно выполняться
$$ (u,w,v)=(w,u,v)=(u,v,w)=0$$
Возьмем для $u,v$ одночлены $\eta_1,\eta_2\in G$, что
$\widehat{\eta_1}=\widehat{u}$, $\widehat{\eta_2}=\widehat{v}$ \ и
\, $\eta_1\eta_2\ne 0$. И подберем $\eta\in G$,
$\widehat{\eta}=\widehat{w}$ \,такой, чтобы $\eta\eta_1\eta_2\ne
0$. Тогда имеем $u^0=\eta_1\otimes u$, \,$v^0=\eta_2\otimes v$,
\,$w^0=\eta\otimes w \in G(J)$ --- элементы грассмановой оболочки,
а $w^0$ лежит в ее центре $G(J)$, значит,
$$ (u^0,w^0,v^0)=(w^0,u^0,v^0)=(u^0,v^0,w^0)=0.$$
т.\,е. получаем
$$0=(\eta\otimes w,\,\eta_1\otimes u,\,\eta_2\otimes v)=
\eta\eta_1\eta_2\otimes(w,u,v),$$ откуда $(w,u,v)=0$. Точно таким
же образом $(u,v,w)=0$. Лемма доказана.


\medskip

{\bf Следствие 1.} {\em Пусть $J$ --- конечномерная простая
унитальная йорданова супералгебра с полупростой четной частью
над алгебраически замкнутым полем $P$, отличная от супералгебры
(\ref{komp v'}). Тогда все тернарные дифференцирования $J$
являются стандартными, т.\,е. вида (\ref{komp2})}.


\medskip

{\bf Доказательство.} Пусть $\Delta=(D,F,G)\in TDer(J)$.
Как было показано в лемме 1,
$F=D-L_g$, \ $G=D-L_f$, \ где $f=F(1)$,
 $g=G(1)$, и в итоге \ $f=c+w/2$, \ $g=c-w/2$, \ где \
$c=D(1)/2$, \ а $w$ удовлетворяет равенству $(x,w,y)=0$ для любых
$x,y\in J$, иначе говоря, $(J,w,J)=0$, откуда по лемме 3 получаем
$w\in Z(J)= P$, т.\,е. $w$
--- скаляр из поля $P$.

Cлучай 1. $\widehat{D}=\widehat{F}=\widehat{G}=0$. В силу теоремы
1, главная компонента $D$ в этом случае явлется стандартной,
т.\,е. $D=\alpha I + D^0$, где $\alpha\in P$ --- скаляр, а $D^0\in
Der(J)$ --- действует по правилу обычного дифференцирования в
обычных алгебрах, и соответственно, $D^0(1)=0$. Значит, в этом
случае имеем $c=\alpha/2$, и соответственно, $f=(\alpha+w)/2 \in
P$, \ $g=(\alpha-w)/2 \in P$
--- скаляры из поля, и тогда
$$F=D-L_g=D-gI=\alpha I + D^0 - \frac{1}{2}(\alpha-w)I = \frac{1}{2}(\alpha+w)I + D^0,$$
$$G=D-L_f=D-fI=\alpha I + D^0 - \frac{1}{2}(\alpha+w)I = \frac{1}{2}(\alpha-w)I + D^0.$$
Обозначив \,$\alpha I=\varphi$, \ $(\alpha+w)I/2=\chi$, \
$(\alpha-w)I/2=\psi$,\, получаем как раз (\ref{komp2}) для
$(D,F,G)$, при этом $\varphi=\chi+\psi$, $D^0(xy)=D^0(x)y +
xD^0(y).$

Cлучай 2. $\widehat{D}=\widehat{F}=\widehat{G}=1$. Тогда, как было
показано в теореме 1, $D(1)=0=2c$ и тогда $F(1)=f=c+w/2=w/2\in P$,
\ $G(1)=g=c-w/2=-w/2\in P$ --- скаляры из поля, т.\,е. лежат в
четной части $J_0$. В то же время при нечетных $F,G$ должно быть
$F(1),G(1)\in J_1$; значит, остается $F(1)=G(1)=0=w=f=g$, и тогда
уже  $F=G=D=D^1$, где $D^1$, в свою очередь, нечетное
дифференцирование, действующее по правилу
$$
D^1(xy)=D^1(x)y+(-1)^{\widehat{x}}xD^1(y).
$$
Следствие доказано.

\medskip

Рассмотрим теперь случай (vi) из основной классификационной
теоремы. Определим трехмерную супералгебру Капланского $K_3$:
$$(K_3)_0=Pe, \ \  (K_3)_1=Pz+Pw$$
где $e^2=e$, $ez=\frac{1}{2}z$, $ew=\frac{1}{2}w$, $zw=e$.

\medskip

{\bf Теорема 2.} {\em Всякое обобщенное дифференцирование $D$ в
cупералгебре $K_3$ является стандартным, т.\,е.  вида
(\ref{komp3}).}


\medskip

{\bf Доказательство.} Случай 1. $D$ --- четное, т.\,е.
$\widehat{D}=0$. Тогда основное равенство (1) \mbox {для $D$:}
$$ D(xy)=F(x)y + xG(y),$$

Рассмотрим образы $d=D(e)$, $f=F(e)$, $g=G(e)$ четного элемента
$e\in (K_3)_0$ --- они лежат снова в четной части $(K_3)_0$, и
значит, $d=\gamma_0 e$, $f=\gamma_1 e$, $g=\gamma_2 e$. При этом
$d=f+g$, т.\,к. $D(e)=D(e^2)=F(e)e+eG(e)=fe+eg=f+g$, и
соответственно, $\gamma_0=\gamma_1+\gamma_2$.

Отметим следующую формулу при умножении произвольного однородного
элемента $x\in K_3$ на четный базисный элемент $e\in (K_3)_0$:
$$ex=xe=\delta(\widehat{x})x,$$
где $\delta(\widehat{x})=1\ \mbox{при}\ \widehat{x}=0$, и\,
$\delta(\widehat{x})=1/2\ \mbox{при}\ \widehat{x}=1$; (например,
взять $\delta(\widehat{x})=1-\frac{1}{2}\widehat{x}$).

Положим в равенстве (1) $y=e$. Тогда
$$ D(xe)=F(x)e+xg;$$
$$\delta(\widehat{x})D(x)=\delta(\widehat{F(x)})+\gamma_2 xe.$$
Т.\,к. $\widehat{F}=0$, то $\widehat{F(x)}=\widehat{x}$, и
соответственно, имеем
$$\delta(\widehat{x})D(x)=\delta(\widehat{x})F(x)+\delta(\widehat{x})\gamma_2x,$$
$$ D(x)=F(x)+\gamma_2 x, \ \ F(x)=D(x)-\gamma_2 x, \ \
F=D-\gamma_2 I.$$ C другой стороны,
$$x=e \ \ \Longrightarrow \ \ D(ey)=fy + eG(y);$$
$$\delta(\widehat{y})D(y)=\gamma_1 ey+\delta(\widehat{G(y)})G(y), \ \
\delta(\widehat{y})D(y)=\gamma_1\delta(\widehat{y})y+\delta(\widehat{y})G(y),$$
$$ D(y)=\gamma_1 y + G(y), \ \ G(y)=D(y)-\gamma_1 y, \ \
G=D-\gamma_1 I.$$

Подставляя  выражения $F,G$ снова в (\ref{komp}), получаем
$$
D(xy) = (D(x)-\gamma_2x)y + x(D(y)-\gamma_1y)=
D(x)y+xD(y)-\gamma_0xy,
$$
и тогда, положив $D^0=D-\gamma_0 I$, имеем
$$D^0(xy)=D(xy)-\gamma_0xy=D(x)y+xD(y)-2\gamma_0xy=D^0(x)y+xD^0(y),$$
т.\,е. $D^0$ --- обычное четное дифференцирование, а
$D=D^0+\gamma_0I$, что и требовалось.

\medskip

Случай 2. $D$ --- нечетное, т.\,е. $\widehat{D}=1$. Тогда основное
равенство (1) \mbox {для $D$:}
$$ D(xy)=F(x)y + (-1)^{\widehat{x}}xG(y)$$

Рассмотрим образы базисных элементов $e\in (K_3)_0$, и $z,w\in
(K_3)_1$, которые в данном случае при отображениях $D,F,G$ должны
менять четность, т.\,е. $D(e),F(e),G(e)\in (K_3)_1$, а
$D(z),F(z),G(z),D(w),F(w),G(w)\in (K_3)_0$, и таким образом,
$$
D(e)=v,\ \ D(z)=\alpha e, \ \ D(w)=\beta e,
$$
$$
F(z)=\gamma_1e, \ \ F(w)=\delta_1e,\ \ G(z)=\gamma_2e, \ \
G(w)=\delta_2e.
$$

Из равенства $e=zw=-wz$ имеем:
\begin{eqnarray}\label{komp1*}
v=D(e)=D(zw)=F(z)w-zG(w)=\frac{1}{2}(\gamma_1w-\delta_2z)=
\end{eqnarray}
$$=-D(wz)=F(w)z-wG(z)=\frac{1}{2}(-\delta_1z+\gamma_2w),$$
откуда $\gamma_1=\gamma_2=\gamma$, \,$\delta_1=\delta_2=\delta$, и
соответственно, $F(z)=G(z)=\gamma e$, $F(w)=G(w)=\delta e$, т.\,е.
$F\equiv G|_{J_1}$.

Далее из равенства $ez=ze=\frac{1}{2}z$ имеем:
\begin{eqnarray}\label{komp2*}
\frac{1}{2}D(z)=D(ez)=F(e)z+eG(z)=F(e)z+\gamma e=
\end{eqnarray}
$$=D(ze)=F(z)e-zG(e)=\gamma e+G(e)z.$$
Отсюда видим, что $(F(e)-G(e))z=0$, что возможно только в случае
$(F(e)-G(e))=k_1z$ (и тогда, поскольку $z\in J_1$
--- нечетное, имеем $z^2=0$).
Аналогично $(F(e)-G(e))w=0$, и $(F(e)-G(e))=k_2w$. Значит,
$F(e)-G(e)=0$, т.\,е. $F(e)=G(e)$.

Наконец, из равенства $e^2=e$ имеем:
$$v=D(e)=D(e^2)=F(e)e+eG(e)=2F(e)e=F(e)$$
(последнее в силу того, что $F(e)\in J_1=L(z,w)$). Таким образом,
$D(e)=F(e)=G(e)=v$.

Из (\ref{komp2*}) имеем
$$\alpha e/2=vz+\gamma e.$$
В то же время, из (\ref{komp1*}) $v=(\gamma w-\delta z)/2$,
\,тогда \,$vz=\gamma w/2=-\gamma e/2$, и подставляя, имеем:
$$\alpha e/2=-\gamma e/2 + \gamma e,$$
откуда $\alpha=\gamma$, и тогда $D(z)=F(z)=G(z)=\alpha e$.

Аналогично $\delta=\beta$, и соответственно, $D(w)=F(w)=G(w)=\beta
e$.

Таким образом, действие $D,F,G$ совпадает на $e,z,w$, т.\,е. на
всех базисных элементах, а значит, $D,F,G$ совпадают на всей
супералгебре. Тогда имеем

$$D(xy)=F(x)y+(-1)^{\widehat x}xG(y)=D(x)y+(-1)^{\widehat x}xD(y),$$
и т.о. $D=D^1$ --- обыкновенное нечетное дифференцирование, ч.т.д.
Теорема доказана.

\medskip

\texttt{ЗAМЕЧАНИЕ.} Отметим, что получили $v=(\alpha w-\beta
z)/2$,
и таким образом, мы доказали также, что всякое нечетное обобщенное
дифференцирование $D$ супералгебры $K_3$ имеет вид:
\begin{eqnarray}\label{komp3*}
D(e)=-\frac{1}{2}\beta z+\frac{1}{2}\alpha w, \ \ D(z)=\alpha e, \
\ D(w)=\beta e,
\end{eqnarray}
где $\alpha,\beta$ --- произвольные
скаляры.

\medskip

{\bf Следствие 2.} \emph{Всякое тернарное дифференцирование
$\Delta=(D,F,G)$ супералгебры Капланского $K_3$ является
стандартным, т.\,е. вида (\ref{komp2}).}

\medskip

{\bf Доказательство.} Cлучай 1.
$\widehat{D}=\widehat{F}=\widehat{G}=0$. Как показано в теореме 2,
в этом случае $F=D-\gamma_2 I$, $G=D-\gamma_1 I$, в свою очередь,
главная компонента $D$ явлется стандартной, $D=D^0+\gamma_0 I$,
где $\gamma_1+\gamma_2=\gamma_0$
--- скаляры из $P$, а $D^0\in Der(J)$ --- стандартное четное
дифференцирование. Таким образом,

$$F=D-\gamma_2 I = D^0 + (\gamma_0-\gamma_2) I $$
$$G=D-\gamma_1 I = D^0 + (\gamma_0-\gamma_1) I $$
Обозначив $\gamma_0 I=\varphi$, \, $(\gamma_0-\gamma_2)I=\chi$, \
$(\gamma_0-\gamma_1)I=\psi$,\, получаем как раз (\ref{komp2}) для
$(D,F,G)$, $\varphi=\chi+\psi$.

Cлучай 2. $\widehat{D}=\widehat{F}=\widehat{G}=1$. Тогда, как было
показано в теореме 2, \,$F=G=D=D^1$ --- стандартное нечетное
дифференцирование.
Следствие доказано.

\medskip

Рассмотрим алгебру Грассмана $\Gamma$ от (нечетных)
антикоммутирующих порождающих $\xi_1,\ldots,\xi_n,\ldots$. Для
определения нового умножения используем операцию
$$
\frac{\partial}{\partial
\xi_j}(\xi_{i_1}\xi_{i_2}\cdots\xi_{i_n})= \left\{
\begin{array}{ll}(-1)^{k-1}\xi_{i_1}\cdots
\xi_{i_{k-1}}\xi_{i_{k+1}}\cdots
\xi_{i_n}\ \ \text{при}\ j=i_k, \\
0\ \ \text{при}\ j\neq i_s\ \forall\, s=1,\ldots n
\end{array}
\right.
$$

Для $f,g\in \Gamma_0\cup\Gamma_1$ грассманово умножение
определяется следующим образом:
$$
\{ f,g\} =(-1)^{\widehat{f}}\,\sum_{j=1}^{\infty}\frac{\partial
f}{\partial e_j}\frac{\partial g}{\partial e_j}.
$$

Пусть $\overline{\Gamma}$ --- изоморфная копия $\Gamma$ с
отображением изоморфизма $x \rightarrow \overline{x}$. Рассмотрим
прямую сумму $J(\Gamma)=\Gamma+\overline{\Gamma}$ векторных
пространств и зададим на ней структуру йордановой супералгебры,
положив $A_0= \Gamma_0+\overline\Gamma_1$, $A_1=
\Gamma_1+\overline\Gamma_0$ с умножением $\bullet$ :
$$
a\bullet b = ab, \ \ \overline{a}\bullet b = (-1)^{\widehat b}\,
\overline{ab\,}, \ \ a\bullet\overline{b\,}=\overline{ab\,}, \ \
\overline{a}\overline{b\,}=(-1)^{\widehat b}\,\{ a,b\},
$$
где $a,b\in\Gamma_0\cup\Gamma_1$ и $ab$ --- произведение в
$\Gamma$. Пусть $\Gamma_n$ --- подалгебра алгебры $\Gamma$,
порожденная элементами $\xi_1,\ldots,\xi_n$. Через $J(\Gamma_n)$
обозначим подсупералгебру $\Gamma_n+\overline{\Gamma_n}$
супералгебры $J(\Gamma)$. Если $n\geq 2$, то $J(\Gamma_n)$
--- простая йорданова супералгебра.

\medskip

{\bf Теорема 3.} {\em Все oбобщенные дифференцирования
супералгебры $J(\Gamma_n)$ являются стандартными, т.\,е. вида
(\ref{komp3}).}

\medskip

{\bf Доказательство.} Пусть $\Gamma=F<\xi_1,\ldots,\xi_n\ |\
\xi_i\xi_j=-\xi_j\xi_i>$
--- алгебра Грассмана, $J=\Gamma+\overline{\Gamma}$. Положим
$$e_1=\frac{1+\overline\xi_k}{2},\ e_2=\frac{1-\overline\xi_k}{2}$$
тогда $e_i\in J_0$ --- четные, при этом выполнены соотношения
$$e_i^2=e_i, \ e_ie_j=0, \ e_1+e_2=1,$$
т.\,е. $e_1,e_2$ --- система ортогональных идемпотентов, значит
имеет место пирсовское разложение  $J=J_{11}+J_{12}+J_{22}$.

Для произвольного $a\in J$ представимо в виде
$a=\gamma^1+\xi_k\gamma^2+\overline{\delta^1}+\overline{\xi_k\delta^2}$,
где $\gamma^i, \delta^i\in \Gamma$ и
$\partial\gamma^i/\partial\xi_k=\partial\delta^i/\partial\xi_k=0$,
(т.\,е. все $\gamma^i,\delta^i$ не содержат $\xi_k$). Имеем

$$ a\cdot\overline\xi_k=\overline\xi_k\gamma^1+\delta^2$$
$$ (a\cdot\overline\xi_k)\cdot\overline\xi_k=\gamma^1+\overline\xi_k\delta^2$$

Отметим, что

$$
a\in J_{11}\ \ \Longleftrightarrow\ \
ae_1=a\ \ \Longleftrightarrow\ \ a\cdot\overline\xi_k=a\ \
\Longleftrightarrow\ \
a=(1+\overline\xi_k)\gamma,\ \ \partial\gamma/\partial\xi_k=0;
$$
$$
a\in J_{12}\ \ \Longleftrightarrow\ \
ae_1=\frac{1}{2}a\ \ \Longleftrightarrow\ \
a\cdot\overline\xi_k=0\ \
\Longleftrightarrow\ \
a=\xi_k\gamma+\overline{\delta},\ \
\partial\gamma/\partial\xi_k=\partial\delta/\partial\xi_k=0;
$$
$$
a\in J_{22}\ \ \Longleftrightarrow\ \
ae_2=a\ \ \Longleftrightarrow\ \ a\cdot\overline\xi_k=-a\ \
\Longleftrightarrow\ \
a=(1-\overline\xi_k)\gamma,\ \ \partial\gamma/\partial\xi_k=0.
$$

Пусть теперь $D$ --- произвольное обобщенное дифференцирование
$J$. По лемме 1, которая в данном случае также в силе,
$$D(e_1)=D(e_1^2)=2D(e_1)e_1 -2(e_1c)e_1, \ \ 2c=D(1)$$
таким образом, $D(e_1)\in J\cdot e_1\subseteq J_{11}\oplus
J_{12}$, соответственно, $D(e_1)=a_{11}+a_{12}$; \,аналогично
$D(e_2)=a_{22}+a_{21}$. Тогда
$$2c=D(1)=D(e_1)+D(e_2)=a_{11}+a_{22}+a_{12}+a_{21}.$$
Возвращаясь, имеем
$$
D(e_1)= 2D(e_1)e_1 -2e_1c\cdot e_1=
2a_{11}+a_{12}-((a_{11}+a_{22}+a_{12}+a_{21})\cdot e_1)\cdot e_1 =
$$
$$
2a_{11}+a_{12}-(a_{11}+\frac{1}{2}(a_{12}+a_{21}))\cdot e_1 =
a_{11}+a_{12}-\frac{1}{4}(a_{12}+a_{21}),
$$
значит, $a_{12}+a_{21}=0$.

Заметим, что \,$\overline\xi_k\cdot e_1=e_1$, \
$\overline\xi_k\cdot J_{12}=0$, \ $\overline\xi_k\cdot e_2=-e_2$,
\,и тогда
$$
(\xi_k,e_1,a_{12})=e_1\cdot
a_{12}-\overline\xi_k\cdot(e_1a_{12})=\frac{1}{2}(a_{12}-\overline\xi_ka_{12})=\frac{1}{2}a_{12},\
\ \ (\xi_k,e_2,a_{12})=\frac{1}{2}a_{21}.
$$
Положим $D'=D-2D_{\overline\xi_k,a_{12}}$, тогда имеем
$D'(e_1)=a_{11}$, $D'(e_2)=a_{22}$. Значит, теперь без ограничения
общности можем считать $a_{12}=a_{21}=0$.

Итак, имеем $2c=D(1)=a_{11}+a_{22}$. Заметим еще, что
$\overline\xi_k=e_1-e_2$, и таким образом, кроме того,
$D(\overline\xi_k)=a_{11}-a_{22}$.

$$
D(\overline\xi_k)=D(\xi_k\cdot\overline{1})=
D(\xi_k)\cdot\overline{1}+(-1)^{\widehat{D}}\xi_k\cdot
D(\overline{1}) - (-1)^{\widehat{D}}(\xi_k c \cdot\overline{1} +
\xi_k\cdot c\overline{1})
$$

Пусть теперь, как уже отмечалось,
$a_{11}=(1+\overline\xi_k)\gamma_1$,
$a_{22}=(1+\overline\xi_k)\gamma_2$, где $\gamma_i\in\Gamma$ и
$\partial\gamma_i/\partial\xi_k=0$.
Тогда имеем
$$2c=a_{11}+a_{22}=(\gamma_1+\gamma_2)+\overline\xi_k(\gamma_1-\gamma_2),$$
$$
2\xi_kc\cdot\overline{1}=(\xi_k(\gamma_1+\gamma_2))\cdot\overline{1}=\overline{\xi_k\gamma_1}+\overline{\xi_k\gamma_2},
\ \ 2\xi_k\cdot c\overline{1}=
\xi_k(\overline{\gamma_1}+\overline{\gamma_2})=\overline{\xi_k\gamma_1}+\overline{\xi_k\gamma_2}.
$$
Возвращаясь и подставляя, получаем
$$
D(\overline\xi_k)=D(\xi_k\cdot\overline{1})=
D(\xi_k)\cdot\overline{1}+(-1)^{\widehat{D}}\xi_k\cdot
D(\overline{1}) -
(-1)^{\widehat{D}}\overline{\xi_k(\gamma_1+\gamma_2)}.
$$
С другой стороны,
$$D(\overline\xi_k)=a_{11}-a_{22}=(\gamma_1-\gamma_2)+\overline\xi_k(\gamma_1+\gamma_2).$$
Отметим, что четность $\widehat{D} = \widehat{D(1)} = \widehat{c}
= \widehat{\gamma_1} = \widehat{\gamma_2}$, и соответственно,

$$
(-1)^{\widehat{D}}(\overline{\xi_k(\gamma_1+\gamma_2)})=(-1)^{\widehat{\gamma_1}}(\overline{\xi_k(\gamma_1+\gamma_2)})=\overline\xi_k(\gamma_1+\gamma_2).
$$
Значит,
\begin{eqnarray}\label{komp10}
(\gamma_1-\gamma_2)+2\overline\xi_k(\gamma_1+\gamma_2)=
D(\xi_k)\cdot\overline{1}+(-1)^{\widehat{D}}\xi_k\cdot
D(\overline{1}).
\end{eqnarray}

Заметим, что $\overline\xi_k(\gamma_1+\gamma_2),\,
D(\xi_k)\cdot\overline{1}\in \overline{\Gamma}$. Далее $\xi_k\cdot
D(\overline{1})\in\xi_k\cdot(\Gamma+\overline{\Gamma})$, поэтому
$$\xi_k\cdot D(\overline{1})=\xi_k\cdot\alpha+\overline{\beta},$$
где $\alpha,\beta\in\Gamma$, $\partial\alpha/\partial\xi_k=0$.

Сравнивая обе части (\ref{komp10}) видим, что
$\gamma_1-\gamma_2=(-1)^{\widehat{D}}\cdot\xi_k\alpha$,\, откуда
$\alpha=0$, $\gamma_1=\gamma_2=\gamma$.

Итак, $c=\gamma$, т.\,е. $\partial c/\partial\xi_k=0$ --- это
значит, $c\in P$, $D(1)=2c\in P$, значит, $D$ --- стандартное.

\medskip

{\bf Следствие 3.} \emph{Все тернарные дифференцирования
супералгебры $J(\Gamma_n)$ являются стандартными, т.\,е. вида
(\ref{komp2})}.

\medskip

{\bf Доказательство.} Аналогично следствию 1.

\medskip

Подытоживая, можно сформулировать следующие результаты.

\medskip

{\bf Теорема 4.} {\em Пусть $J=J_0\oplus J_1$,
$J_1\ne 0$ --- конечномерная простая йорданова супералгебра,
отличная от супералгебры (\ref{komp v'}), над алгебраически
замкнутым полем характеристики 0. Тогда все oбобщенные
(соответственно, тернарные) дифференцирования в $J$ являются
стандартными, т.\,е. вида (\ref{komp3}) (соответственно,
(\ref{komp2}))}.


\medskip

{\bf Доказательство.} Следует из теорем [21]-[22],1,2,3 и
следствий 1,2,3.

\medskip

В случае произвольной характеристики не равной 2, описание
обобщенных и тернарных дифференцирований получено при
дополнительном ограничении полупростоты четной части:

\medskip

{\bf Теорема 5.} {\em Пусть $J=J_0\oplus J_1$, $J_1\ne 0$ ---
конечномерная простая йорданова супералгебра с полупростой четной
частью $J_0$, отличная от супералгебры (\ref{komp v'}), над
алгебраически замкнутым полем характеристики не равной 2. Тогда
все oбобщенные (соответственно, тернарные) дифференцирования $J$
являются стандартными, т.\,е. вида (\ref{komp3}) (соответственно,
(\ref{komp2}))}.

\medskip

{\bf Доказательство.} Следует из теорем [20],1,2 и следствий 1,2.



\medskip


Рассмотрим теперь ситуацию, исключенную до сих пор из
расссмотрения
--- супералгебру (\ref{komp v'}) билинейной суперформы $J(V,f)$ с
одномерной четной частью:
$$J=Pe+V,\ \ V=V_0+V_1, \ \ V_0=0,$$
cоответственно тогда $J_0=Pe$, $J_1=V_1=V$.

\medskip

{\bf Теорема 6.} {\em Пусть
$J(V,f)=Pe+V$
--- супералгебра невырожденной суперформы
с одномерной четной частью $J_0=Pe$ над алгебраически замкнутым
полем $P$. Тогда все четные обобщенные дифференцирования в $J$
являются стандартными
т.\,е. вида (\ref{komp3}).
Нечетные же ненулевые обобщенные дифференцирования в $J$
существуют только при размерности векторной части $dim\,V=2$, ---
в этом случае все они являются нестандартными и описываются
следующим образом:}
$$
\{\,D_v: J\longrightarrow J\ |\ D_v(1)=v,\
D_v(x)=\frac{1}{2}f(x,v)\ \,\forall\,x\in V \,\}_{v\in V}.
$$

\medskip

{\bf Доказательство.} Пусть $D$ --- произвольное обобщенное
дифференцирование

Случай 1. $D$ --- четное, т.\,е. $\widehat{D}=0$. Тогда $D(1)\in
J_0=Pe$, \,т.\,е. $D(1)=\alpha\cdot 1$, $\alpha\in P$,  откуда
сразу следует, что $D$ --- стандартное четное дифференицрование.


\medskip

Случай 2. $D$ --- нечетное, т.\,е. $\widehat{D}=1$. В силу
кососимметричности $f$ на $V$ существует базис
$u_1,v_1,\ldots,u_n,v_n$ пространства $V$, такой, что
$f(u_i,v_i)=1,\ f(u_i,v_j)=0\ \forall\,i\ne j$. Т.\,к.
$\widehat{D}=1$, то $D(u_i),D(v_i)\in J_0=Pe$, т.\,е.
$D(u_i)=\alpha_i e$, $D(v_i)=\beta_i e$, \ $\alpha_i,\beta_i\in
P$.

Применим лемму 1:
$$D(1)=D(u_iv_i)=\alpha_iv_i-\beta_iu_i+(f(u_i,c)v_i +
f(c,v_i)u_i)=\gamma_iu_i+\delta_iv_i\in L(u_i,v_i)$$

Если $n>1$, т.\,е. $dim\,V>2$, то
$\exists\,\{u_j,v_j\}\neq\{u_i,v_i\}$, и тогда
$D(1)=\gamma_ju_j+\delta_jv_j\in L(u_j,v_j)$. Но поскольку
$u_i,v_i,u_j,v_j$ --- базисные вектора и линейно независимы, то
такое возможно только при $D(1)=0$. Соответственно, $c=D(1)/2=0$ и
тогда
$$
0=D(1)=D(u_iv_i)=\alpha_iv_i-\beta_iu_i+((u_ic)v_i +
u_i(cv_i))=\alpha_iv_i-\beta_iu_i,
$$
откуда в силу линейной независимости $u_i,v_i$ имеем
$\alpha_i=\beta_i=0$, и таким образом, $D\equiv 0$ на $J$, т.\,е.
при размерности $dim\, V>2$ в $J$ нет ненулевых нечетных
дифференцирований.

Осталось рассмотреть случай $n=1$, т.\,е. $dim\,V=2$.

Пусть $D$
--- произвольное нечетное обобщенное дифференцирование.
Допустим, $D(1)=v\in V$. Тогда $\exists\, u\in V:\ uv=f(u,v)=1$.
Запишем
$$D(1)=D(uv)=D(u)v-uD(v)+\frac{1}{2}((uD(1))v+u(D(1)v));$$
$$v=D(u)v-uD(v)+\frac{1}{2}(f(u,v)v+uf(v,v)),$$
$$v=-D(v)u+(D(u)+1/2)v,$$
откуда $D(v)=0$, $D(u)=1/2$.

Рассмотрим ограничение $D|_V=l: V\longrightarrow P$ --- это
линейный функционал значит, в силу невырожденности формы $f$ на
$V$
$$\exists\, v_0\in V:\ l(x)=f(x,v_0).$$
Тогда $f(v,v_0)=D(v)=0$ --- отсюда, опять же, в силу
невырожденности $f$ на $V=L(u,v)$ следует, что $v_0=kv$. Теперь
имеем: $1/2=D(u)=f(u,v_0)=kf(u,v)=k.$ Итак, получаем в итоге
$$D(x)=l(x)=f(x,v_0)=kf(x,v)=\frac{1}{2}f(x,v)\ \ \forall\,x\in V.$$

Можно провести и обратно: пусть $v\in V$ --- произвольно.
Определим отображение $D: J\longrightarrow J$ по правилу $D(1)=v,\
D(x)=f(x,v)/2\ \,\forall\,x\in V$. И проверим равенство (\ref{komp4}) 
для $D$:
$$D(xy) =
D(x)y+(-1)^{\widehat{x}}xD(y)-\frac{(-1)^{\widehat{x}}}{2}((xD(1))y+x(D(1)y))\
\ \ \forall \,x,y \in J.
$$
При $x\in Pe$ или $y\in Pe$ --- очевидно, выполняется. При $x,y\in
V$ данное равенство принимает вид:
$$f(x,y)v=f(x,v)y-f(y,v)x.$$
Поскольку $dim\,V=2$, то $v=\alpha x+\beta y$, или $x=\lambda y$,
после чего нетрудно убедиться, что равенство выполнено в силу
кососимметричности $f$. Значит, по лемме 1, $D$ --- обобщенное
дифференцирование.

Таким образом, можно установить взаимно однозначное соответствие
между всеми векторами $v\in V$ и всеми нечетными обобщенными
дифференцированиями $D$ по правилу
$$D(1)=v\,;\ D(x)=\frac{1}{2}f(x,v)\ \ \forall\, x\in V.$$
Теорема доказана.

\medskip

{\bf Следствие 4.} {\em Пусть
$J(V,f)=Pe+V$
--- супералгебра невырожденной суперформы
с одномерной четной частью $J_0=Pe$ над алгебраически замкнутым
полем $P$. Тогда все четные тернарные дифференцирования в $J$
являются стандартными
т.\,е. вида (\ref{komp2}).
Нечетные же ненулевые тернарные дифференцирования в $J$ существуют
только при размерности векторной части $dim\,V=2$, --- в этом
случае все они являются нестандартными и описываются следующим
образом:}
$$
\{\,\Delta_v\ |\ \Delta_v(1)=(v,\frac{1}{2}v,\frac{1}{2}v),\
\Delta_v(x)=(\frac{1}{2}f(x,v),f(x,v),f(x,v))\,\ \forall\,x\in V
\,\}_{v\in V}.
$$

\medskip

{\bf Доказательство.} Пусть $\Delta=(D,F,G)$ --- произвольное
тернарное дифференцирование $J$.

Cлучай 1. $\widehat{D}=\widehat{F}=\widehat{G}=0$. См. аналогичный
случай следствия 1.

Случай 2. $\widehat{D}=\widehat{F}=\widehat{G}=1$. Всё точно так
же, как в следствии 1, за исключением равенства $D(1)=2c=0$, --- в
нашем случае $D(1)=v$ из теоремы 6, также имеем $D(x)=f(x,v)/2\ \
\forall\,x\in V$. В остальном, здесь так же имеют место
соотношения: $F=D-L_g$, \ $G=D-L_f$, \ где $f=F(1)$, $g=G(1)$, и в
итоге $f=c+w/2$, $g=c-w/2$,  где  $c=D(1)/2$, а
$w$ --- скаляр из поля $P$. Таким образом, $f-g=w\in P$; в то же
время при нечетных $F,G$ должно быть $f-g=(F-G)(1)\in J_1$.
Значит, остается $f-g=0$, и $f=g=c=D(1)/2=v/2$. Тогда
$$F(x)=G(x)=D(x)-gx=\frac{1}{2}f(x,v)-\frac{1}{2}vx=f(x,v).$$
Следствие доказано.

\medskip

В заключение, обобщим теорему 4 на полупростые супералгебры.
Используем следующий результат:

\medskip

{\bf Теорема \cite{Zelmanov2}.} {\em Пусть $J$ --- конечномерная
йорданова супералгебра над полем $P$ характерстики не 2.
Супералгебра $J$ полупроста тогда и только тогда, когда $J$ есть
прямая сумма простых йордановых супералгебр и унитальных блоков
$H_K(J_1\oplus\ldots\oplus J_s)=J_1\oplus\ldots\oplus J_s + K\cdot
1$, где $J_1,\ldots, J_s$ --- простые неунитальные йордановы
супералгебры над некоторым расширением $K$ поля $P$, иными
словами,
$$
J=\bigoplus_{i=1}^s (J_{i1}\oplus\ldots\oplus J_{ir_i} + K_i\cdot
1)\oplus J_1\oplus\ldots\oplus J_t,
$$
где $J_1,\ldots, J_t$ --- простые йордановы $P$-супералгебры,
$K_1,\ldots, K_s$ --- расширения поля $P$, и для каждого $i$
$J_{i1},\ldots, J_{ir_i}$ --- простые неунитальные супералгебры
над $K_i$}.

\medskip

Отметим, что в рассматриваемых нами условиях (характеристика поля
равна 0, или четная часть полупроста) есть только один тип
неунитальных простых супералгебр --- Капланского. Поэтому докажем
сначала следующую

\medskip

{\bf Лемма 4.} Пусть $J=\bigoplus_{i=1}^n K_3^{(i)} + P\cdot 1$
--- унитальное замыкание $n$ супералгебр Капланского $K_3^{(i)}\cong
K_3$. Тогда всякое обобщенное дифференцирование $J$ является
стандартным.

\medskip

{\bf Доказательство.} Пусть $e_i,z_i,w_i$ --- канонический базис
супералгебры $K_3^{(i)}$.

Случай 1. $D$ --- четное, $\widehat{D}=0$. Имеем
$$ J_0=Pe_1\oplus\ldots\oplus Pe_n+P\cdot 1=Pe_1\oplus\ldots\oplus
Pe_n\oplus Pe_0,$$ где $e_0=1-e_1-\ldots-e_n$. Ясно, что
$D(e_i)=\alpha_ie_i$, \,$i=0,1,\ldots,n$, поскольку
$Pe_i=(Pe_i)^2\lhd J_0$. Тогда и
$D(1)=\sum_{i=0}^n\alpha_ie_i=2c$.

Рассмотрим
$$D(z_i)=2D(e_iz_i)=2D(e_i)z_i+2e_iD(z_i)-2((e_ic)z_i+e_i(cz_i)).$$
Поскольку $\widehat{D}=0$ и соответственно, $\widehat
{D(z_i)}=\widehat{z_i}=1$, значит, $D(z_i)\in L(z_i,w_i)$; тогда
заметим, что $D(z_i)e_i=D(z_i)/2$; и еще заметим\,
$z_ie_0=z_i-z_ie_i=z_i/2$. Теперь имеем
$$D(z_i)=D(z_i)+\alpha_iz_i-(\frac{1}{4}(\alpha_i+\alpha_0)z_i+\frac{1}{2}\alpha_iz_i),$$
откуда $\alpha_i=\alpha_0=\alpha$. Тогда $D^0=D-\alpha I$ ---
обыкновенное четное дифференцирование.

Случай 2. $D$ --- нечетное, $\widehat{D}=1$. Поскольку
$K_3^{(i)}=(K_3^{(i)})^2\lhd J$, то нетрудно понять, что
$D(K_3^{(i)})=D((K_3^{(i)})^2)\subseteq K_3^{(i)}\
\,\forall\,i=1,\ldots,n$. В силу замечания после теоремы 3, имеет
место вид (\ref{komp3*}) для $D$, т.\,е.
$$D(e_i)=\beta z_i+\alpha w_i, \ \ D(z_i)=2\alpha_ie_i, \ \ D(w)=-2\beta_ie_i$$
при некоторых скалярах $\alpha_i,\beta_i\in P$,  $i=1,\ldots,n$.

Т.\,к. $\widehat{D(1)}=\widehat{D}=1$, соответственно, $D(1)\in
J_1$, то имеем
$D(1)=\sum_{i=1}^n\lambda_iz_i+\sum_{i=1}^n\mu_iw_i$.

Рассмотрим
$$2\alpha_ie_i=D(z_i)=2D(e_iz_i)=2D(e_i)z_i+2e_iD(z_i)-2[(e_ic)z_i+e_i(cz_i)]=$$
$$=2\alpha_iw_iz_i+4\alpha_ie_i-(e_i(\lambda_iz_i+\mu_iw_i)\cdot z_i
+ e_i\cdot\mu_iw_iz_i)=
2\alpha_ie_i-(-\frac{1}{2}\mu_ie_i-\mu_ie_i),$$ откуда $\mu_i=0$.
Точно так же, аналогично $\lambda_i=0$, и для всех $i=1,\ldots,n$.
Таким образом получаем $D(1)=0$, что сразу означает, что $D$ ---
обыкновенное дифференцирование. Лемма доказана.

\medskip

{\bf Теорема 7.} {\em Пусть $J$ - конечномерная полупростая
йорданова супералгебра над алгебраически замкнутым полем
характеристики 0, $J = J'\oplus J''$,  где в $J''$ собраны все
прямые слагаемые вида (\ref{komp v'}), а $J'$ таких слагаемых не
содержит. Пусть $J''=J(V_1,f_1)\oplus\cdots J(V_k,f_k)$. Тогда
всякое четное обобщенное(/тернарное) дифференцирование
супералгебры $J$ стандартно, а всякое нечетное
обобщенное(/тернарное) дифференцирование $D$ (/$\Delta$) имеет вид
\,$D=D'+D_{v_1}+\cdots+D_{v_k}$
(/$\Delta=\Delta'+\Delta_{v_1}+\cdots+\Delta_{v_k}$), где
$D'$(/$\Delta'$) - нечетное дифференцирование супералгебры $J'$, а
$D_{v_i}$(/$\Delta_{v_i}$) - обобщенное(/тернарное)
дифференцирование супералгебры $J(V_i,f_i)$, $v_i\in V_i$,
$i=1,\ldots,k$. В частности, если $J=J'$, т.\,е. $J$ не содержит
прямых слагаемых вида (\ref{komp v'}), то все
обобщенные(/тернарные) дифференцирования супералгебры $J$
стандартны.}

\medskip


{\bf Доказательство.} Теперь имеем
$$J=J'\oplus J''= (J_1\oplus\ldots\oplus J_s\oplus J_{s+1}\oplus\ldots\oplus J_t)\oplus(J_{t+1}\oplus\ldots\oplus J_{t+k}),$$
где  $J_1,\ldots,J_s$ --- рассматривавшиеся унитальные замыкания
неунитальных алгебр,  $J_{s+1},\ldots,J_t$ --- простые унитальные
супералгебры, за исключением вида (\ref{komp v'}),
$J_{t+1}=J(V_1,f_1),\ldots, J_{t+k}=J(V_k,f_k)$.

Как бы то ни было, все подалгебры $J_i$ являются унитальными
и потому $J_i=J_i^2\lhd J$ и соответственно, $D(J_i)\subseteq
D(J_i^2) \subseteq J_i$,
--- действует инвариантно на $J_i$.
Таким образом, $D_i=D|_{J_i}$ является обобщенным
дифференцированием во всех подалгебрах $J_i$. Тогда для
$J_1,\ldots,J_s$ по лемме 4, для $J_{s+1},\ldots,J_t$ --- по
теореме 4, а для $J_{t+1},\ldots,J_{t+k}$--- по теореме 6,
действие $D_i$ в четном случае является стандартным. Значит, имеем
$$
D=\sum_{i=1}^{t+k} D_i = \sum_{i=1}^{t+k} (\alpha_iI_i + D_i^0)=
\sum_{i=1}^{t+k} \alpha_ie_iI + \sum_{i=1}^{t+k} D_i^0=\alpha I +
D^0,
$$
где $\alpha_i\in P$ --- скаляры из поля, $I_i=I|_{J_i}$, $D_i^0\in
Der(J_i)$
--- обыкновенные \text{дифференцирования в $J_i$}, $e_i$ --- единицы в
$J_i$, и $\alpha=\sum_{i=1}^{t+k} \alpha_ie_i$. Тогда очевидно,
$\alpha\in Z(J)$
--- элемент из центра $J$, а $D^0=\sum_{i=1}^{t+k} D_i^0\in
Der(J)$
--- обыкновенное дифференцирование всей супералгебры, и таким образом, $D$
--- стандартное четное обобщенное дифференцирование.

В нечетном случае --- аналогично действие $D'=D|_{J'}=\sum_{i=1}^t
D_i$ является стандартным, т.\,е. обыкновенным нечетным
дифференцированнием, а для $i=t+1,\ldots,t+k$ имеем $D_i=D_{v_i}$
из теоремы 6, и таким образом, получаем искомый вид
$D=D'+D_{v_1}+\cdots+D_{v_k}$.

Для тернарных дифференцирований --- доказательство следует из
только что доказанного для обобщенных, аналогично следствию 1.

\bigskip

ШЕСТАКОВ Алексей Иванович\\
Институт математики им. С.Л.Соболева СО РАН, пр. Академика Коптюга, 4;\\
Новосибирск, 630090, РОССИЯ.\\
e-mail: shestalexy@yandex.ru, \\
тел. (383)363-45-57, (+7)913-753-3738

\end{document}